\documentclass[a4paper, 11pt,reqno]{amsart} 
\usepackage{amssymb}
\usepackage{amsmath}

\newtheorem{thm}{Theorem}[section]

\newtheorem{rem}{Remark}[section]
\newtheorem{cor}{Corollary}[section]

\makeatletter
   
          \@addtoreset{equation}{section}
\makeatother

\begin{document}

\baselineskip=1.15\baselineskip   

\title[Massless Dirac operators]{The asymptotic limits of zero modes 
of   massless Dirac operators }

\author{Yoshimi Sait\={o}
and
 Tomio Umeda\\}

\address{
Department of Mathematics \\
University of Alabama at Birmingham, 
Birmingham, AL 35294, USA}
\email{saito@math.uab.edu}

\address{Department of Mathematical Sciences,
University of Hyogo, Himeji 671-2201, Japan
}
\email{umeda@sci.u-hyogo.ac.jp}
\date{}

\maketitle

\vspace{20pt}

\textbf{Abstract.} Asymptotic behaviors of zero modes   of 
the massless Dirac operator $H=\alpha\cdot  D + Q(x)$ are discussed, where
 $\alpha= (\alpha_1, \, \alpha_2, \, \alpha_3)$ is
the triple of  $4 \times 4$ Dirac matrices,
$ D=\frac{1}{\, i \,} \nabla_x$, and 
$Q(x)=\big( q_{jk} (x) \big)$   is a $4\times 4$ Hermitian matrix-valued function
with
 $| q_{jk}(x) | \le C \langle x \rangle^{-\rho} $,  $\rho >1$.
 We shall show that for  every zero mode $f$, 
the asymptotic limit of  $|x|^2f(x)$
 as $|x| \to +\infty$ exists.
 The limit is expressed in terms of an integral of $Q(x)f(x)$.

 \vspace{15pt}


 \vspace{15pt}
\textbf{The 2000 Mathematical Subject Classification:} 35Q40, 35P99, 81Q10

\newpage

\section{Introduction}

In this paper we study asymptotic behaviors of zero modes 
(i.e., eigenfunctions with
the zero eigenvalue; see Definition 1.1) of 
the massless Dirac operator
\begin{equation} \label{eqn:1-1}
H= \alpha \cdot D + Q(x), \quad D=\frac{1}{\, i \,} \nabla_x,
\,\,\, x \in {\mathbb R}^3,
\end{equation}
where $\alpha= (\alpha_1, \, \alpha_2, \, \alpha_3)$ is
the triple of  $4 \times 4$ Dirac matrices
\begin{equation*}\label{eqn:1-2}
\alpha_j = 
\begin{pmatrix}
 \mathbf 0 &\sigma_j \\ \sigma_j &   \mathbf 0 
\end{pmatrix}  \qquad  (j = 1, \, 2, \, 3)
\end{equation*}
with  the $2\times 2$ zero matrix $\mathbf 0$ and
the triple of  $2 \times 2$ Pauli matrices
\begin{equation*}\label{eqn:1-3}
\sigma_1 =
\begin{pmatrix}
0&1 \\ 1& 0
\end{pmatrix}, \,\,\,
\sigma_2 =
\begin{pmatrix}
0& -i  \\ i&0
\end{pmatrix}, \,\,\,
\sigma_3 =
\begin{pmatrix}
1&0 \\ 0&-1
\end{pmatrix},
\end{equation*}
and $Q(x)$ is a $4\times 4$ Hermitian matrix-valued function
decaying at infinity.

Particular emphasis must be placed on the fact  that
one can view the operator (\ref{eqn:1-1}) as 
a generalization of 
the operator
\begin{equation} \label{eqn:1-4}
\alpha\cdot \big(D - A(x) \big) + q(x) I_4 , 
\end{equation}
where $(q, A)$ is an electromagnetic potential and
$I_4$ is the $4\times 4$ identity matrix,
by  taking $Q(x)$ to be $- \alpha\cdot A(x) + q(x) I_4$.
In the case  where $q(x)\equiv 0$, the operator
(\ref{eqn:1-4}) becomes of the form
\begin{equation}   \label{eqn:1-4-1}
\alpha\cdot \big(D - A(x) \big)
=
\begin{pmatrix}
 \mathbf 0 &\sigma \cdot (D -A(x)) \\ 
\sigma \cdot (D -A(x)) & \mathbf 0 
\end{pmatrix}.
\end{equation}
The component $\sigma \cdot (D -A(x))$ is known as 
the Weyl-Dirac operator. See Balinsky and Evans \cite{BalinEvan2}.

The paper by Fr\"ohlich, Lieb 
and Loss \cite{FrohlichLiebLoss} revealed 
that the existence of zero modes
 of a  
Weyl-Dirac operator
plays  a crucial role 
in the study of stability of Coulomb systems with magnetic
fields.
In connection with  \cite{FrohlichLiebLoss},
Loss and Yau \cite{LossYau} 
constructed, for the first time ever,   
examples  of  vector potentials $A(x)$ for which
 the corresponding Weyl-Dirac 
operators have zero modes. 
After the work by Loss and Yau \cite{LossYau}
was published,
there have been many contributions on the study of
zero modes of Weyl-Dirac operators. 
See Adam, Muratori and Nash \cite{AdamMuratoriNash1},
\cite{AdamMuratoriNash2}, \cite{AdamMuratoriNash3},
Balinsky and Evans \cite{BalinEvan1},
\cite{BalinEvan2},
\cite{BalinEvan3},
Bugliaro, Fefferman and Graf \cite{BugliaroFeffermanGraf},
Elton \cite{Elton}
and,
Erd\"os and Solovej \cite{ErdosSolovej1}, 
 \cite{ErdosSolovej2}, 
\cite{ErdosSolovej3}.

We would like to mention 
Loss and Yau's example of the zero mode $\psi_L$ and 
the vector potential $A_L$:
\begin{align}    
\psi_{L}(x) &= \langle x \rangle^{-3}
 \big( I_2 + i \sigma \cdot x \big) \phi_0 ,   \label{eqn:LY-1}\\
\noalign{\vskip 2pt}
A_L(x) &= 3 \langle x \rangle^{-4} 
  \big\{ (1-|x|^2)  w_0 + 2 (w_0 \cdot x) x + 2 w_0 \times x \big\},  \label{eqn:LY-2}
\end{align}
where
$\langle x \rangle = \sqrt{1 + |x|^2 \,}$,
$\phi_0 = {}^t(1, \, 0)$, 
and
\begin{equation*}
w_0= \big( \phi_0 \cdot (\sigma_1 \phi_0), \, 
  \phi_0 \cdot (\sigma_2 \phi_0), \,
  \phi_0 \cdot (\sigma_3 \phi_0) \big).
\end{equation*}
It follows from (\ref{eqn:LY-1}) that
\begin{equation}   \label{eqn:LY-3}
\lim_{r\to +\infty} r^2 \psi_{L}(r \omega) = (i \sigma \cdot \omega) \phi_0,
\end{equation}
where $r= |x|$ and $\omega =x/ |x|$.

Adam, Muratori and Nash \cite{AdamMuratoriNash1},
\cite{AdamMuratoriNash2}, \cite{AdamMuratoriNash3} 
developed the idea of Loss and Yau \cite{LossYau} and
constructed many examples of 
the pairs of zero modes and vector potentials in a systematic way. 
Among other things, it is important in the context of the present paper 
that they constructed the
zero modes of the form 
\begin{equation}  \label{eqn:amn-1}
\psi_{A}(x) = \langle x \rangle^{-2} U(x) \phi_0,
\end{equation}
where $U(x)$ is a  $2 \times 2$ matrix-valued
function with the limit
\begin{equation}   \label{eqn:amn-2}
U_{\infty}(\omega):= \lim_{r \to +\infty} U(r\omega).
\end{equation}
Thus, it follows from (\ref{eqn:amn-1}) and (\ref{eqn:amn-2})
that $r^2\psi_{A}(r\omega)$ has a limit as $r\to +\infty$:
\begin{equation}  \label{eqn:amn-3}
\lim_{r\to + \infty}r^2\psi_{A}(r\omega)
 =  U_{\infty}(\omega) \phi_0.
\end{equation}

It is apparent from (\ref{eqn:LY-3})
and (\ref{eqn:amn-3}) that both 
$\psi_L(x)$ and $\psi_A(x)$ behave
 in the same manner as  $r\to +\infty$.
We would like to emphasize that
this is not a sheer  coincidence.
Actually, Theorem \ref{thm:th-lim} below
asserts
that every zero mode $\psi(x)$ of the Weyl-Dirac operator
behaves like
\begin{equation}  \label{eqn:insert-1}
\psi(r\omega) \sim r^{-2} \, i(\sigma \cdot \omega) \psi_0
\;\;\; \quad 
(\psi_0 \in {\mathbb C^2} \mbox{\ \rm a constant vector})
\end{equation}
for $r\to +\infty$
if the vector potential $A$ satisfies
$|A(x)|\le const. \langle x \rangle^{-\rho} \;\;\;(\rho >1)$.

The purpose of the present paper
is to show that every zero mode $f(x)$ 
of the massless Dirac operator (\ref{eqn:1-1})
behaves like
\begin{equation}  \label{eqn:insert-2}
f(r\omega) \sim r^{-2} \, i(\alpha \cdot \omega) f_0
\;\;\; \quad 
(f_0 \in {\mathbb C^4} \mbox{\ \rm a constant vector})
\end{equation}
for $r\to +\infty$
if each component of $Q(x)$ satisfies
the inequality (\ref{eqn:2-1}) in Assumption (A) below;
see Theorem
\ref{thm:th-lim}. 
We should like to note that 
Theorem \ref{thm:th-lim} can be regarded
as a refinement of
our previous result \cite[Theorem 2.1]{SaitoUmeda}, where
we proved that every 
zero mode $f(x)$ 
of the operator (\ref{eqn:1-1}) 
satisfies the inequality 
\begin{equation}    \label{eqn:previous}
|f(x)| \le const. \langle x \rangle^{-2}
\end{equation}
under the same assumption as in the 
present paper.
\vspace{15pt}

\noindent
\textbf{Notation.}

\noindent
By $L^2 =L^2({\mathbb R}^3)$, we mean the Hilbert space of
square-integrable functions on ${\mathbb R}^3$, and 
we introduce  a Hilbert space ${\mathcal L}^2$ by
     ${\mathcal L}^2 = [L^2({\mathbb R}^3)]^4$, where 
the inner product  
is given by
\begin{equation*}
    (f, g)_{{\mathcal L}^2}
 = \sum_{j=1}^4 (f_j, g_j)_{L^2}
\end{equation*}
for  $f = {}^t(f_1, f_2, f_3, f_4)$ 
and
$g = {}^t(g_1, g_2, g_3, g_4)$.
By $H^{1}({\mathbb R}^3)$ we denote
the Sobolev space of order $1$,
and by  
${\mathcal H}^{1}$ we mean the
Hilbert space $[H^{1}({\mathbb R}^3)]^4$.
When we mention the Weyl-Dirac operator, we must handle
 two-vectors (two
components spinors)
which will be denoted by $\psi$.

\vspace{15pt}

\noindent
\textbf{Assumption (A).} 

\noindent
Each element $q_{jk}(x)$ 
($j, \, k =1, \, \cdots, \, 4$) of $Q(x)$ is 
a measurable function satisfying
\begin{equation} \label{eqn:2-1}
| q_{jk}(x) | \le C_{\! q} \langle x \rangle^{-\rho} 
\quad   ( \rho >1 ),
\end{equation}
where $C_{\! q}$ is a positive constant.
Moreover, $Q(x)$ is a Hermitian matrix for each 
$x \in {\mathbb R^3}$.

\vspace{15pt}

Note that, under Assumption (A), the Dirac operator (\ref{eqn:1-1})
is a self-adjoint operator in $\mathcal L^2$ with  
$\mbox{Dom}(H) = {\mathcal H}^1$.
The self-adjoint realization will be denoted 
by $H$ again.

\vspace{15pt}
\noindent
D{\scriptsize EFINITION} 1.1. \ 
By a zero mode, we mean a function 
$f \in \mbox{Dom}(H)$ which satisfies
$Hf=0$.

\vspace{15pt}

We are now in a position to state the main result
of the present paper.

\vspace{10pt}

\begin{thm} \label{thm:th-lim}
Suppose Assumption {\rm(A)} is satisfied. Let $f$ be  
a zero mode of the operator {\rm(\ref{eqn:1-1})}.
Then for any $\omega \in {\mathbb S}^2$ 
\begin{equation}   \label{eqn:lim-1}
\lim_{r\to +\infty}r^{2}f(r\omega)
= 
-\frac{i}{\, 4\pi \,} \,
(\alpha \cdot \omega) \! \!
\int_{{\mathbb R}^3}   Q(y) f(y) \, dy,
\end{equation}
where the convergence being uniform
with respect to $\omega \in {\mathbb S}^2$.
\end{thm}

\vspace{10pt}

In connection with the expression $f(r\omega)$ in (\ref{eqn:lim-1}), 
it is worthy to
note that every zero mode is a continuous function
(see Theorem \ref{thm:th-SUmain} in the beginning of section 2).

Since $\alpha \cdot \omega$
 is a unitary matrix, we have an immediate corollary to 
Theorem \ref{thm:th-lim}.

\begin{cor}  \label{cor:unitary}
For any $\omega \in {\mathbb S}^2$
\begin{equation}   \label{eqn:lim-1-1}
\lim_{r\to +\infty}r^{2} |f(r\omega)|
= 
\frac{1}{\, 4\pi \,} \,
\Big|
\int_{{\mathbb R}^3}   Q(y) f(y) \, dy
\Big|.
\end{equation}
\end{cor}

\vspace{10pt}

One should note  that Corollary \ref{cor:unitary} assures the
$\omega$-independence of the limit of $r^{2}
|f(r\omega)|$  for $r \to \infty$.
In particular we see that
Corollary \ref{cor:unitary} implies an interesting fact:
\begin{equation}   \label{eqn:lim-1-2}
\lim_{r\to +\infty}r^{2}f(r\omega)=0 \;
\mbox{ for some (any) }\omega
\Longleftrightarrow
\int_{{\mathbb R}^3}  Q(y) f(y) \, dy =0.
\end{equation}

\vspace{10pt}
As for a zero mode of the Weyl-Dirac
operator, we have
the following theorem,
which is also a corollary to Theorem \ref{thm:th-lim}.

\begin{thm} \label{thm:th-limWD}
Suppose 
\begin{equation} \label{eqn:2-1WD}
| A(x) | \le C \langle x \rangle^{-\rho} 
\quad   ( \rho >1 ),
\end{equation}
where $C$ is a positive constant.
Let $\psi$ be  
a zero mode of the Weyl-Dirac
operator $\sigma\cdot(D-A(x))$.
Then for any $\omega \in {\mathbb S}^2$ 
\begin{gather}
\begin{split}   \label{eqn:lim-1WD}
{}&\lim_{r\to +\infty}r^{2}\psi(r\omega)     \\
&= 
\frac{i}{\, 4\pi \,}
\int_{{\mathbb R}^3}  
\big\{ 
  \big( \omega\cdot A(y) \big)  I_2
    + i\sigma \cdot 
    \big( \omega\times A(y) \big) 
   \big\}  \psi(y) \, dy,
\end{split}
\end{gather}
where the convergence being uniform
with respect to $\omega \in {\mathbb S}^2$.
\end{thm}

Erd\"os and Solovey \cite{ErdosSolovej1} generalized 
the examples by Loss and Yau  \cite{LossYau} from the 
  geometrical point of view, and proposed an intrinsic way
of producing magnetic fields on
${\mathbb S}^3$ and ${\mathbb R}^3$ for which the 
corresponding Weyl-Dirac operators have
zero modes. They did not mention asymptotic properties of
their zero modes, which were obviously not their concern though. 

It is interesting from our point of view that
Elton \cite{Elton} showed that for 
any integer $m \ge 0$ and an open subset 
$\Omega \subset {\mathbb R}^3$
there exists a vector potential 
$A \in [C_0^{\infty}({\mathbb R}^3)]^3$ such
that $\mbox{supp }A \subset \Omega$ and 
the corresponding Weyl-Dirac operator  has 
a degeneracy of zero modes with multiplicity $m$.
This fact, together with Theorem \ref{thm:th-limWD}, 
indicates that the asymptotic behavior of 
vector potential $A$ does not affect the asymptotic behavior of
zero modes of the corresponding Weyl-Dirac operator
as long as $A$ satisfies the hypothesis (\ref{eqn:2-1WD}).

\section{Proofs}

The proof of Theorem \ref{thm:th-lim} is 
based on an estimate,
 which was established 
in our previous paper \cite[Theorem 2.1]{SaitoUmeda}.

\vspace{10pt}

\begin{thm}[Sait\={o} and Umeda] \label{thm:th-SUmain}
Suppose Assumption {\rm(A)} is satisfied. Let $f$ be  
a zero mode of the operator {\rm(\ref{eqn:1-1})}.
Then 

{\rm (i)} the  inequality 
\begin{equation}   \label{eqn:2-2}
|f(x)| \le C
\langle x \rangle^{-2}
\end{equation}
holds for all $x \in {\mathbb R}^3$, where the constant 
$C(=C_f)$ depends only on the zero mode $f$;

\vspace{4pt}

{\rm (ii)} the zero mode $f$ is a continuous function on ${\mathbb R}^3$.
\end{thm}
\vspace{15pt}

Also, the proof of Theorem \ref{thm:th-lim} is 
based on a fact that every zero mode $f$ of 
the operator {\rm(\ref{eqn:1-1})} satisfies 
the integral equation
\begin{equation} \label{eqn:sio-1-lim}
f(x) = 
 - \displaystyle{
\frac{i}{\,4\pi\,}
\int_{{\mathbb R}^3}
   \frac{\alpha \cdot (x - y)}{ |x-y|^3} 
\, Q(y) f(y) \, dy.
}
\end{equation}
This fact was established in our previous paper \cite{SaitoUmeda} too; 
see
(5.3) in Section 5 of \cite{SaitoUmeda}.

\vspace{5pt}

\begin{rem}
If  we formally take the limit of \hbox{\rm(\ref{eqn:sio-2-lim})} below 
as $r \to +\infty$, then we can
readily obtain \hbox{\rm(\ref{eqn:lim-1})}. 
Unfortunately, this argument is not rigorous.
\end{rem}

\vspace{5pt}

\begin{proof}[{\bf Proof of Theorem \ref{thm:th-lim}}]
We begin with the integral equation (\ref{eqn:sio-1-lim}) with
$x=r\omega\;\; (\omega \in {\mathbb S}^2)$, and multiply the both
sides of (\ref{eqn:sio-1-lim}) by $r^2$:
\begin{equation} \label{eqn:sio-2-lim}
r^2 f(r\omega) = 
 - \displaystyle{
\frac{i}{\,4\pi\,}
\int_{{\mathbb R}^3}
   \frac{\alpha \cdot (\omega - r^{-1}y)}{|\omega-r^{-1}y|^3} 
\, Q(y) f(y) \, dy.
}
\end{equation}
We then see from (\ref{eqn:sio-2-lim}) that
\begin{gather}   \label{eqn:sio-3-lim}
\begin{split}
{}&r^2 f(r\omega) 
+
\frac{i}{\,4\pi\,}
\int_{{\mathbb R}^3}
  (\alpha \cdot \omega )
\, Q(y) f(y) \, dy     \\
&\;=
\frac{i}{\,4\pi\,}
\int_{{\mathbb R}^3}
\alpha \cdot
\Big\{ \omega -
   \frac{\omega - r^{-1}y}{\,|\omega-r^{-1}y|^3} \Big\}
\, Q(y) f(y) \, dy.
\end{split}
\end{gather}

Now let $\varepsilon >0$ be given so that $0< \varepsilon <1/2$, 
and choose
$R_0$ so that
\begin{equation}    \label{eqn:epsilon-1}
R_0^{-\rho + 1 }  < \varepsilon.
\end{equation}
Note that $\rho > 1$; see Assumption (A).
For $r \ge 2R_0$, we define
\begin{align}
E_1 &:=\big\{ \, y \in {\mathbb R}^3 \, 
    \big| \;    |y| \le R_0 \, \big\}, 
\label{eqn:set-1}  \\ 
E_2 &:=  \big\{ \, y \in {\mathbb R}^3 \, 
    \big| \;   |y| >  R_0, \,  |r\omega - y| \le \frac{r}{2} \, \big\}, 
\label{eqn:set-2} \\ 
E_3 &:=   \big\{ \, y \in {\mathbb R}^3 \, 
    \big| \;   |y| >  R_0, \,  |r\omega - y| > \frac{r}{2} \, \big\},  
\label{eqn:set-3} 
\end{align}
and accordingly we decompose the integral on the right hand side
of (\ref{eqn:sio-3-lim}) into three parts:
\begin{align}  
I_r(\omega)& := 
\frac{i}{\,4\pi\,}
\int_{E_1}
\alpha \cdot
\Big\{ \omega -
   \frac{\omega - r^{-1}y}{\,|\omega-r^{-1}y|^3} \Big\}
\, Q(y) f(y) \, dy,   
\label{eqn:integral-1}  \\
I\!I_r(\omega)& := 
\frac{i}{\,4\pi\,}
\int_{E_2}
\alpha \cdot
\Big\{ \omega -
   \frac{\omega - r^{-1}y}{\,|\omega-r^{-1}y|^3} \Big\}
\, Q(y) f(y) \, dy,   
\label{eqn:integral-2}  \\
I\!I\!I_r(\omega)& := 
\frac{i}{\,4\pi\,}
\int_{E_3}
\alpha \cdot
\Big\{ \omega -
   \frac{\omega - r^{-1}y}{\,|\omega-r^{-1}y|^3} \Big\}
\, Q(y) f(y) \, dy.   
\label{eqn:integral-3}  
\end{align}
We thus have
\begin{gather}
\begin{split}   \label{eqn:decomposition-1}
{}&r^2 f(r\omega) 
+
\frac{i}{\,4\pi\,}
\int_{{\mathbb R}^3}
  (\alpha \cdot \omega )
\, Q(y) f(y) \, dy   \\
\noalign{\vskip 4pt}
{}&= I_r(\omega) + I\!I_r(\omega) + I\!I\!I_r(\omega).
\end{split}
\end{gather}

To estimate $I_r(\omega)$, we first note that
\begin{equation}   \label{eqn:decomposition-2}
\frac{1}{\, 2 \,} \le 
 |\omega-r^{-1}y|  \le
    \frac{3}{\, 2 \,} 
\quad \mbox{if } |y| \le R_0, \; r\ge 2R_0.
\end{equation}
This implies that
\begin{gather}
\begin{split}  \label{eqn:decomposition-3}  
\Big| \omega -
   \frac{\omega - r^{-1}y}{\,|\omega-r^{-1}y|^3} \Big|
&\le
\Big( \frac{1}{\, 2 \,} \Big)^{\!-3}
\big|\,
|\omega-r^{-1}y|^3 \omega - (\omega-r^{-1}y) 
 \,\big|   \\
\noalign{\vskip 4pt}
&=
2^3 \,
\big| \,
(|\omega-r^{-1}y|^3  - 1 )\omega  + r^{-1}y 
\,\big|
\end{split}
\end{gather}
when $|y| \le R_0, \; r\ge 2R_0$.
Moreover, we have
\begin{gather}
\begin{split}  \label{eqn:decomposition-4}  
\big| \,
|\omega-r^{-1}y|^3  - 1 
\,\big|
&=
\big| \,
|\omega-r^{-1}y|  - 1 
\,\big|     \\
& \qquad  \times
\big( |\omega-r^{-1}y|^2 + |\omega-r^{-1}y| + 1   \big)   \\
\noalign{\vskip 4pt}
&\le
\big| \,
|\omega-r^{-1}y|  - 1 
\,\big| 
\times \frac{\, 19 \,}{4}   \\
\noalign{\vskip 4pt}
&=
\frac{\, 19 \,}{4} 
\frac{1}{\, |\omega-r^{-1}y|  + 1  \,}
\times 
\big| \,
|\omega-r^{-1}y|^2  - 1 
\,\big|    \\
\noalign{\vskip 4pt}
&\le
\frac{\, 19 \,}{6} 
\big| \,
-2\omega \cdot (r^{-1}y) +r^{-2}|y|^2  
\,\big|     
\end{split}
\end{gather}
provided that $|y| \le R_0, \; r\ge 2R_0$.
Combining (\ref{eqn:decomposition-3}) and (\ref{eqn:decomposition-4}),
we obtain
\begin{align} \label{eqn:decomposition-5}
\Big| \omega -
   \frac{\omega - r^{-1}y}{\,|\omega-r^{-1}y|^3} \Big|
\le 
2^3
\Big\{
\frac{\, 19 \,}{6} (2R_0 r^{-1} + R_0^2 r^{-2} ) 
+ R_0 r^{-1}
\Big\}
\end{align}
whenever $|y| \le R_0, \; r\ge 2R_0$.
Now it follows from (\ref{eqn:decomposition-5}),
 Theorem \ref{thm:th-SUmain}, Assumption (A) 
and the anti-commutation relation
that
\begin{gather}
\begin{split} \label{eqn:decomposition-6}
| I_r(\omega) |
&\le 
\frac{1}{\, 4\pi \,}
\times
2^3
\Big\{
\frac{\, 22 \,}{3} R_0 r^{-1} + \frac{\, 19 \,}{6}R_0^2 r^{-2} 
\Big\}
\int_{E_1}
|Q(y) f(y)| \, dy  \\
&\le
C_1  R_0 r^{-1}
\int_{{\mathbb R}^3} 
  \langle y \rangle^{-\rho -2} \, dy   \\
\noalign{\vskip 2pt}
&\le
C_1^{\prime}  R_0 r^{-1}
\end{split}
\end{gather}
for all $r \ge 2R_0$ and all $\omega \in {\mathbb S}^2$,
where the constant $C_1^{\prime}$ is
dependent only on the constant $C_f$ in Theorem \ref{thm:th-SUmain}
and the constant $C_{\! q}$ in Assumption (A).
Note that in the first inequality in (\ref{eqn:decomposition-6}) 
we have used the fact that $|(\alpha \cdot x)f| = |x| \, |f|$ for
all $x \in {\mathbb R}^3$ and all $f \in {\mathbb C}^4$,
and that in the third inequality we have
used the fact that 
$\langle y \rangle^{-\rho -2}$ is integrable on ${\mathbb R}^3$
since $\rho +2 > 3$.

As for $I\!I_r(\omega)$, it follows again from
Theorem \ref{thm:th-SUmain} and Assumption (A) 
that
\begin{equation}   \label{eqn:decomposition-7}
| I\!I_r(\omega) |
\le 
C_2
\int_{E_2} 
\Big( 1 +
   \frac{1}{ \, |\omega-r^{-1}y|^2 \,} \Big)
\, \langle y \rangle^{-\rho -2} \, dy
\end{equation}
for all $r\ge 2R_0$ and all $\omega \in {\mathbb S}^2$,
where the constant $C_2$ 
depends only on the constants $C_f$ and $C_{\! q}$.
To estimate the right hand side of (\ref{eqn:decomposition-7}),
we need the fact that
\begin{equation}    \label{eqn:decomposition-8}
y \in E_2 \, \Rightarrow  \; |y| \ge \frac{r}{\, 2 \,}.
\end{equation}
Thus, the integral on the right hand side of (\ref{eqn:decomposition-7})
is estimated by
\begin{gather}   
\begin{split}    \label{eqn:decomposition-9}
{}& C \Big\{
\int_{|y|\ge r/2}
  \langle y \rangle^{-\rho -2} \, dy 
+
r^{-\rho -2}  
  \int_{|r\omega - y|\le r/2} 
  \frac{r^2}{ \, |r\omega- y|^2 \,} \, dy 
\Big\}    \\
\noalign{\vskip 4pt}
{}& \quad
\le
C^{\prime} r^{-\rho +1}
\end{split}
\end{gather}
for for all $r\ge 2R_0$ and all $\omega \in {\mathbb S}^2$,
with the constant $C^{\prime}$ 
 independent of $\omega$ and $r$.
Combining (\ref{eqn:decomposition-9}) with
(\ref{eqn:decomposition-7}), we get
\begin{equation}     \label{eqn:decomposition-10}
| I\!I_r(\omega) |
\le 
C_2^{\prime}  \,  r^{-\rho +1}
\end{equation}
for for all $r\ge 2R_0$ and all $\omega \in {\mathbb S}^2$,
where  $C_2^{\prime}$ is
a constant independent of $\omega$ and $r$.

In the same way as in (\ref{eqn:decomposition-7}) we have
\begin{equation}   \label{eqn:decomposition-11}
| I\!I\!I_r(\omega) |
\le 
C_3
\int_{E_3} 
\Big( 1 +
   \frac{1}{ \, |\omega-r^{-1}y|^2 \,} \Big)
\, \langle y \rangle^{-\rho -2} \, dy
\end{equation}
for all $r\ge 2R_0$ and all $\omega \in {\mathbb S}^2$.
Here the constant $C_3$ 
depends only on  $C_f$  $C_{\! q}$.
The integral on the right hand side of (\ref{eqn:decomposition-11})
is bounded by
\begin{gather}   
\begin{split}    \label{eqn:decomposition-12}
{}& C^{\prime\prime} \Big\{
\int_{|y|\ge R_0}
  \langle y \rangle^{-\rho -2} \, dy 
+
  \int_{E_3} 
  \frac{1}{ \, |\omega- r^{-1}y|^2 \,} 
     \langle y \rangle^{-\rho -2} \, dy 
\Big\}    \\
\noalign{\vskip 4pt}
{}& \quad
\le
C^{\prime\prime\prime} R_0^{-\rho +1}
\end{split}
\end{gather}
for for all $r\ge 2R_0$ and all $\omega \in {\mathbb S}^2$,
with the constant $C^{\prime\prime\prime}$ 
 independent of $\omega$ and $r$.
Here we have used the fact that $|\omega- r^{-1}y|\ge 1/2$
for all $y \in E_3$. 
It follows from (\ref{eqn:decomposition-11}) and 
(\ref{eqn:decomposition-12}) that
\begin{equation}   \label{eqn:decomposition-13}
| I\!I\!I_r(\omega) |
\le 
C_3^{\prime}  R_0^{-\rho +1}
\end{equation}
for all $r\ge 2R_0$ and all $\omega \in {\mathbb S}^2$,
with the constant $C_3^{\prime}$ 
 independent of $\omega$ and $r$.\break

We are now ready to combine (\ref{eqn:decomposition-1})
with (\ref{eqn:decomposition-6}), (\ref{eqn:decomposition-10})
 (\ref{eqn:decomposition-13}), and we can conclude that
\begin{gather}
\begin{split}   \label{eqn:decomposition-14}
\Big|
r^2 f(r\omega) 
+
\frac{i}{\,4\pi\,}
\int_{{\mathbb R}^3}
  (\alpha \cdot \omega )
\, Q(y) f(y) \, dy 
\Big|    \\
\noalign{\vskip 3pt}
\le 
C_1^{\prime}  R_0 r^{-1} + C_2^{\prime}  \,  r^{-\rho +1}
+C_3^{\prime}  R_0^{-\rho +1}   \\
\le 
C_1^{\prime}  R_0 r^{-1} +
( C_2^{\prime} + C_3^{\prime} )R_0^{-\rho +1} 
\end{split}
\end{gather}
for all $r\ge 2R_0$ and all $\omega \in {\mathbb S}^2$.
Putting $R_1:=R_0/\varepsilon(>2R_0)$, and recalling
(\ref{eqn:epsilon-1}), we have  shown that
\begin{equation}   \label{eqn:umedaconclusion-1}
\Big|
r^2 f(r\omega) 
+
\frac{i}{\,4\pi\,}
\int_{{\mathbb R}^3}
  (\alpha \cdot \omega )
\, Q(y) f(y) \, dy 
\Big|  
\le 
(C_1^{\prime}  + C_2^{\prime} + C_3^{\prime} ) \varepsilon
\end{equation}
for for all $r\ge R_1$ and all $\omega \in {\mathbb S}^2$.
Since $\varepsilon > 0$ was arbitrary, 
(\ref{eqn:umedaconclusion-1}) implies the conclusion
of the theorem.
\end{proof}

\vspace{5pt}

\begin{proof}[{\bf Proof of Theorem \ref{thm:th-limWD}}]
In view of (\ref{eqn:1-4-1}) and Theorem \ref{thm:th-lim}, 
we only have to 
compute
\begin{equation}  \label{eqn:limWD-p1}
-\frac{i}{\, 4\pi \,}
\int_{{\mathbb R}^3}  (\sigma \cdot \omega) 
\big( -\sigma\cdot A(y) \big) \psi(y) \, dy.
\end{equation}
Using
the 
anti-commutation relation 
$\sigma_j \sigma_k  + \sigma_k \sigma_j = 2 \delta_{jk} I_2$
and the facts that
\begin{equation}  \label{eqn:limWD-p2}
\sigma_1\sigma_2= i \sigma_3,
\;\;
\sigma_2\sigma_3= i \sigma_1,
\;\;
\sigma_3\sigma_1= i \sigma_2,
\end{equation}
we get
\begin{equation*}    \label{eqn:limWD-p3}
(\sigma \cdot \omega) 
\big( \sigma\cdot A(y) \big)
= \big( \omega \cdot A(y) \big) I_2
 + i \sigma \cdot
 \big( \omega \times A(y) \big).
\end{equation*}
This completes the proof.
\end{proof}

%%%%%%%%%%%%%%%%%%%%%%%%%%%%%%%%%%%%%%%%%%%%%%%%%%%%%%%%%
%%%%%%%%%%%%%%%%%%%%%%%%%%%%%%%%%%%%%%%%%%%%%%%%%%%%%%%%%
%%%%%%%%%%%%%%%%%%%%%%%%%%%%%%%%%%%%%%%%%%%%%%%%%%%%%%%%%

\end{document}